\documentclass{article}
\usepackage{amsmath} 
\usepackage{amsfonts}
\usepackage{amsthm}
\usepackage{amssymb}
\usepackage{amscd}
\usepackage{tikz}
\usepackage{tikz-cd}
\newtheorem{Thm}{Theorem}

\newtheorem{Lem}{Lemma}
\newtheorem{Rem}{Remark}
\newtheorem{Cor}{Corollary}

\title{Homological Lagrangian monodromy and Wang exact sequence}
\author{Yoshihiro Sugimoto}
\date{}

\begin{document}

\maketitle

\begin{abstract}
    In this paper, we study homological monodromy of a Lagrangian submanifold. We show that homological Lagrangian monodromy is trivial if Hofer norm of a Hamiltonian isotopy is smaller than the minimum energy of J-holomorphic spheres and discs. 
\end{abstract}

\section{Introduction}
Let ${L\subset M}$ be a closed Lagrangian submanifold of a closed or tame symplectic manifold ${(M,\omega)}$. We can consider Hamiltonian diffeomorphisms which preserve $L$ as a set. Such Hamiltonian diffeomorphims form a subgroup of the Hamiltonian diffeomorphism group ${\mathrm{Ham}(M,\omega)}$. We denote this subgroup by ${\mathrm{Ham}_L(M,\omega)}$. An element of ${\mathrm{Ham}_L(M,\omega)}$ naturally induces an automorphism of the singular homology:

\begin{gather*}
    \mathrm{Ham}_L(M,\omega)\longrightarrow \mathrm{Aut}(H_*(L:\mathcal{R})) \\
    \phi \longmapsto \big\{\phi_*:H_*(L:\mathcal{R})\rightarrow H_*(L:\mathcal{R})\big\}.
\end{gather*}
Here, $\mathcal{R}$ is a coefficient ring (usually, ${\mathcal{R}=\mathbb{Z}}$ or ${\mathcal{R}=\mathbb{Z}_2}$). The image of this map ${\mathcal{G}(L:\mathcal{R})}$ is called homological Lagrangian monodromy group. There are several previous results on the properties of ${\mathcal{G}(L:\mathcal{R})}$ (\cite{C2,ASW,HLL,K,MW,O,Y,B}). 

In order to state our main results rigorously, we introduce some notations. Hofer energy of a Hamiltonian path ${\{\phi_t\}\subset \mathrm{Ham}(M,\omega)}$ is defined as follows:
\begin{gather*}
    ||\{\phi_t\}||=\int_0^1\max_{x\in M}H(t,x)-\min_{x\in M}H(t,x)dt \\
    \frac{d}{dt}\phi_t(x)=X_{H_t}(\phi_t(x)).
\end{gather*}
Here, ${X_{H_t}}$ is the Hamiltonian vector field generated by a normalized Hamiltonian function ${H_t=H(t,\cdot)}$. We also define Hofer energy of a homotopy class of path ${\gamma\in \pi_1(\mathrm{Ham}(M,\omega),\mathrm{Ham}_L(M,\omega))}$ as follows:

\begin{gather*}
    ||\gamma||=\inf_{[\{\phi_t\}]=\gamma}||\{\phi_t\}||.
\end{gather*}

Let ${\mathcal{J}(M,\omega)}$ be the set of compatible almost complex structures on ${(M,\omega)}$. We define the minimum energy of holomorphic spheres and discs as follows:
\begin{gather*}
    \sigma_S(M,\omega,J)=\inf\Bigg\{ \ \ \int_{S^2}u^*\omega>0 \ \ \Bigg| \ \ \begin{matrix}u:S^2\rightarrow M \\ \overline{\partial}_Ju=0 \end{matrix}\ \ \Bigg\} 
\end{gather*}
\begin{gather*}
    \sigma_L(M,\omega,J)=\inf\Bigg\{ \ \ \int_{D^2}u^*\omega>0 \ \ \Bigg| \ \ \begin{matrix}u:(D^2,\partial D^2)\rightarrow (M,L) \\ \overline{\partial}_Ju=0 \end{matrix} \ \ \Bigg\}.
\end{gather*}

Any ${\gamma \in \pi_1(\mathrm{Ham}(M,\omega),\mathrm{Ham}_L(M,\omega))}$ determines a symplectic fibration over a disc ${\pi:P\rightarrow D^2}$ and a Lagrangian fibration over the boundary ${\pi:Q\rightarrow \partial D^2}$. See the next section for details of $P$ and $Q$. Let ${TP^{\pi}\subset TP}$ be the kernel of the projection:

\begin{gather*}
    TP^{\pi}=\{v\in TP \ | \ \pi_*(v)=0 \}.
\end{gather*}
Then, ${TP^{\pi}\rightarrow P}$ is a symplectic vector bundle with a symplectic form ${\bar{\omega}}$. We denote the set of almost complex structures on ${TP^{\pi}}$ by ${\mathcal{J}_{\gamma}}$:

\begin{gather*}
    \mathcal{J}_{\gamma}=\Bigg\{\widetilde{J}\in \mathrm{End}(T^{\pi}P) \ \Bigg| \ \begin{matrix}
        \widetilde{J}^2=-\mathrm{Id}  \\ \widetilde{J}:\mathrm{compatible \ with \ } \bar{\omega}
    \end{matrix} \Bigg\}.
\end{gather*}
For any ${\widetilde{J}\in \mathcal{J}_{\gamma}}$, we define ${\sigma_S(\gamma,\widetilde{J})}$ and ${\sigma_L(\gamma,\widetilde{J})}$ as follows:

\begin{gather*}
    \sigma_S(\gamma,\widetilde{J})=\inf \Bigg\{\int_{S^2}u^*\bar{\omega}>0 \ \Bigg| \ \begin{matrix}
        u:S^2\rightarrow P \\ \pi \circ u=\mathrm{constant}, \overline{\partial}_{\widetilde{J}}u=0 
    \end{matrix} \Bigg\}
\end{gather*}
\begin{gather*}
    \sigma_L(\gamma,\widetilde{J})=\inf \Bigg\{\int_{D^2}u^*\bar{\omega}>0 \ \Bigg| \ \begin{matrix}
        u:(D^2,\partial D^2)\rightarrow (P,Q) \\ \pi \circ u=\mathrm{constant}, \overline{\partial}_{\widetilde{J}}u=0 
    \end{matrix} \Bigg\}.
\end{gather*}
In other words, ${\sigma_S(\gamma,\widetilde{J})}$ and ${\sigma_L(\gamma,\widetilde{J})}$ are the minimum energy of ${\widetilde{J}}$-holomorphic spheres and discs in fibers of $P$. Let ${\sigma(M,L,\gamma)}$ be the supremum of them:
\begin{gather*}
    \sigma(M,L,\gamma)=\sup_{\widetilde{J}\in \mathcal{J}_{\gamma}}\min \{\sigma_S(\gamma,\widetilde{J}),\sigma_L(\gamma,\widetilde{J})\}.
\end{gather*}

The next theorem is our main result.
\vspace{7mm}
\begin{Thm}
    Let ${L\subset M}$ be a closed Lagrangian submanifold. Let ${\{\phi_t\}\subset \mathrm{Ham}(M,\omega)}$ be a Hamiltonian isotopy such that ${\phi_0=\mathrm{Id}}$ and ${\phi_1\in \mathrm{Ham}_L(M,\omega)}$ hold. If 
    \begin{gather*}
        ||[\{\phi_t\}]||<\sigma(M,L,[\{\phi_t\}])
    \end{gather*}
    holds, then
    \begin{gather*}
        (\phi_1)_*:H_*(L:\mathbb{Z}_2)\longrightarrow H_*(L:\mathbb{Z}_2)
    \end{gather*}
    is equal to the identity.
\end{Thm}
\vspace{7mm}
A Lagrangian submanifold ${L\subset M}$ is called rational if there is a positive constant ${\eta>0}$ so that ${\omega(\pi_2(M,L))=\eta\mathbb{Z}}$ holds. In this case, ${\sigma(M,L,\gamma)\ge \eta}$ holds for any $\gamma$. So, we have the following corollary.
\vspace{7mm}
\begin{Cor}
    Let ${L\subset M}$ be a closed rational Lagrangian submanifold and let ${\eta>0}$ be the generator of ${\omega(\pi_2(M,L))}$. Let ${\{\phi_t\}\subset \mathrm{Ham}(M,\omega)}$ be a Hamiltonian isotopy such that ${\phi_0=\mathrm{Id}}$, ${\phi_1\in \mathrm{Ham}_L(M,\omega)}$ and ${||\{\phi_t\}||<\eta}$ hold. Then
    \begin{gather*}
        (\phi_1)_*:H_*(L:\mathbb{Z}_2)\longrightarrow H_*(L:\mathbb{Z}_2)
    \end{gather*}
    is equal to the identity.
\end{Cor}

\vspace{7mm}

\begin{Rem}
    In this paper, we use various moduli spaces of $J$-holomoprhic curves. We assume that the transversality is achieved after small perturbations of equations. Today, such perturbation technique is common knowledge among experts. For any moduli space $\mathcal{M}$, ${\sharp_{\mathbb{Z}_2}\mathcal{M}}$ stands for the mod $2$ number of elements of the moduli space in virtial dimension $0$ strata. 
\end{Rem}

\section{Symplectic fibration over $D^2$}
 We fix a Hamiltonian isotopy ${\mathbf{g}=\{g_t\}_{t\in [0,1]}}$ such that ${g_0=\mathrm{Id}}$ and ${g_1\in \mathrm{Ham}_L(M,\omega)}$ hold. Let ${D^2\in \mathbb{C}}$ be the closed unit disc. We divide $D^2$ into positive part ${D^2_+}$ and negative part ${D^2_-}$.
 \begin{gather*}
     D^2_+=\{z\in D^2 \ | \ \mathrm{Re}(z)\ge 0 \} \\
     D^2_-=\{z\in D^2 \ | \ \mathrm{Re}(z)\le 0 \}
 \end{gather*}
The Hamiltonian isotopy ${\mathbf{g}}$ define a symplectic fibration ${\pi:P\rightarrow D^2}$ as follows:

\begin{gather*}
    P=D^2_-\times M\coprod D^2_+\times M/\sim  \\
    (\sqrt{-1}t,g_{\frac{1}{2}(t+1)}(x))\sim (\sqrt{-1}t,x).
\end{gather*}
We also divide the boundary ${\partial D^2}$ into positive part ${\partial^+D^2}$ and negative part ${\partial^-D^2}$.

\begin{gather*}
    \partial^+D^2=\{z\in D^2_+ \ | \ |z|=1\} \\
    \partial^-D^2=\{z\in D^2_- \ | \ |z|=1\}
\end{gather*}
The condition ${g_1\in \mathrm{Ham}_L(M,\omega)}$ implies that there is a Lagrangian fibration ${\pi:Q\rightarrow \partial D^2}$ as follows:

\begin{gather*}
    Q=\partial^-D^2\times L\coprod \partial^+D^2\times L/\sim  \\
    (\sqrt{-1},g_1(x))\sim (\sqrt{-1},x).
\end{gather*}
Wang exact sequence implies that there is a homology long exact sequence of the total space and the fiber.
\begin{gather*}
  \cdots  \longrightarrow H_*(L)\stackrel{(g_1)_*-\mathrm{Id}}{\longrightarrow} H_*(L)\longrightarrow H_*(Q)\longrightarrow H_{*-1}(L)\longrightarrow \cdots
\end{gather*}
So, ${(g_1)_*=\mathrm{Id}}$ holds if and only if the inclusion
\begin{gather*}
    i:H_*(L)\longrightarrow H_*(Q)
\end{gather*}
is injective. For the latter purpose, we introduce a Morse theoretical description of Wang exace sequence. See \cite{BH,S} for a detailed description of Morse homology.

Let ${f:L\rightarrow \mathbb{R}}$ be a Morse function on $L$. We denote the set of critical points of $f$ by ${\mathrm{Crit}(f)}$. Let ${\mathrm{Crit}_k(f)}$ be the set of critical points whose Morse indices are equal to $k$. We fix a Riemannian metric $\rho$ on $L$. We denote the gradient vector field of $f$ with Riemannian metric $\rho$ by ${\mathrm{grad}_{\rho}f}$. The stable manifold and the unstable manifold of critical points ${x\in \mathrm{Crit}(f)}$ are defined as follows:

\begin{gather*}
    W^s(x)=\Bigg\{z\in L \ \ \  \Bigg| \ \ \  \begin{matrix}
     \exists \gamma:[0,\infty)\rightarrow L  \\
     \dot{\gamma}(t)=-\mathrm{grad}_{\rho}f(\gamma(r))  \\
     \gamma(0)=z, \ \lim_{t\rightarrow \infty}\gamma(t)=x
    \end{matrix} \Bigg\} 
\end{gather*}
\begin{gather*}
    W^u(x)=\Bigg\{z\in L \ \ \  \Bigg| \ \ \  \begin{matrix}
     \exists \gamma:(-\infty,0]\rightarrow L  \\
     \dot{\gamma}(t)=-\mathrm{grad}_{\rho}f(\gamma(r))  \\
     \gamma(0)=z, \ \lim_{t\rightarrow -\infty}\gamma(t)=x
    \end{matrix} \Bigg\}.
\end{gather*}
Morse-Smale condition requires that ${W^u(x)}$ and ${W^s(y)}$ intersects transversally for any ${x,y\in \mathrm{Crit}(f)}$. If ${(f,\rho)}$ satisfies Morse-Smale condition, ${(f,\rho)}$ is called a Morse-Smale pair. We assume that ${(f,\rho)}$ is a Morse-Smale pair. Morse chain complex is a $\mathbb{Z}_2$-module generated by ${\mathrm{Crit}(f)}$.

\begin{gather*}
    MC(f,\rho)=\mathbb{Z}_2\langle \mathrm{Crit}(f)\rangle
\end{gather*}
\begin{gather*}
    MC_*(f,\rho)=\mathbb{Z}_2\langle \mathrm{Crit}_{*}(f)\rangle
\end{gather*}
The intersection points ${W^u(x)\cap W^s(y)}$ are canonically identified with the following solution space:

\begin{gather*}
\widetilde{\mathcal{M}}(x,y:(f,\rho))=\Bigg\{\gamma:\mathbb{R}\rightarrow L \ \ \  \Bigg| \ \ \  \begin{matrix}\dot{\gamma}(t)=-\mathrm{grad}_{\rho}f \\ \lim_{t\rightarrow -\infty}\gamma(t)=x, \lim_{t\rightarrow \infty}\gamma(t)=y \end{matrix}\Bigg\}.
\end{gather*}
There is a natural $\mathbb{R}$-action on ${\widetilde{M}(x,y:(f,\rho))}$. If ${(f,\rho)}$ is a Morse-Smale pair,
\begin{gather*}
    \mathcal{M}(x,y:(f,\rho))=\widetilde{\mathcal{M}}(x,y:(f,\rho))/\mathbb{R}
\end{gather*}
is a ${\mathrm{ind}_f(x)-\mathrm{ind}_f(y)-1}$ dimensional manifold. Morse boundary operator ${\partial^M}$ is defined as follows:
\begin{gather*}
    \partial^M:MC_*(f,\rho)\longrightarrow MC_{*-1}(f,\rho)  \\
  \hspace{1cm}  x\longmapsto \sum_{y\in \mathrm{Crit}_{*-1}(f)}\sharp_{\mathbb{Z}_2} \mathcal{M}(x,y:(f,\rho))\cdot y .
\end{gather*}
Morse homology is the homology group of this chain complex.

\begin{gather*}
    MH_*(f,\rho)=H_*(MC(f,\rho),\partial^M)
\end{gather*}

Morse homology does not depend on the various choices in the following meaning. Let ${(f_0,\rho_0)}$ and ${(f_1,\rho_1)}$ be two Morse-Smale pairs on $L$. Let ${\{(l_s,h_s)\}_{s\in \mathbb{R}}}$ be a family of pairs of functions ${l_s\in C^{\infty}(L,\mathbb{R})}$ and Riemannian metrics ${h_s}$ such that

\begin{gather*}
    (l_s,h_s)=\begin{cases}
        (f_0,\rho_0)  & (x\le -R)  \\
        (f_1,\rho_1)  & (x\ge R)
    \end{cases}
\end{gather*}
holds for some ${R>0}$. Let ${\mathcal{N}(x,y:(l_s,h_s))}$ be the following solution space for ${x\in \mathrm{Crit}(f_0)}$ and ${y\in \mathrm{Crit}(f_1)}$.

\begin{gather*}
    \mathcal{N}(x,y:(l_s,h_s))=\Bigg\{ \gamma:\mathbb{R}\rightarrow L \ \ \ \Bigg| \ \ \ \begin{matrix}
        \dot{\gamma}(t)=-\mathrm{grad}_{h_t}l_t(\gamma(t))  \\
        \lim_{t\rightarrow -\infty}\gamma(t)=x, \lim_{t\rightarrow \infty}\gamma(t)=y
    \end{matrix}\Bigg\}
\end{gather*}
For suitable choice of ${(l_s,h_s)}$, ${\mathcal{N}(x,y:(l_s,h_s))}$ is a ${\mathrm{ind}_{f_0}(x)-\mathrm{ind}_{f_1}(y)}$ dimensional manifold for any ${x\in \mathrm{Crit}(f_0)}$ and ${y\in \mathrm{Crit}(f_1)}$. The canonical continuation map between ${MC(f_0,\rho_0)}$ and ${MC(f_1,\rho_1)}$ is defined by counting the number of ${\mathcal{N}(x,y:(l_s,h_s))}$ as follows:

\begin{gather*}
    \Phi_{(f_0,\rho_0),(f_1,\rho_1)}:MC_*(f_0,\rho_0)\longrightarrow MC_*(f_1,\rho_1) \\
    x\longmapsto \sum_{y\in \mathrm{Crit}(f_1)}\sharp_{\mathbb{Z}_2} \mathcal{N}(x,y:(l_s,h_s))\cdot y .
\end{gather*}

Next, we consider Morse thoery of the fibration ${\pi:Q\rightarrow \partial D^2}$. Let ${L_{\pm}=\pi^{-}(\pm 1)}$ be fibers of $\pi$. We choose Morse-Smale pairs ${(f_{\pm},\rho_{\pm})}$ on $L_{\pm}$. We assume that
\begin{gather*}
    \max_{x\in L_-}f_-(x)< \min_{y\in L_+}f_+(y)
\end{gather*}
holds. Let ${F:\rightarrow \mathbb{R}}$ be a Morse function such that 
\begin{gather*}
F|_{L_{\pm}}=f_{\pm}   \\
    \mathrm{Crit}(F)=\mathrm{Crit}(f_-)\coprod \mathrm{Crit}(f_+)
\end{gather*}
holds. We also choose a Riemannian metric ${\widetilde{\rho}}$ on $Q$ so that the following conditions hold.

\begin{itemize}
    \item ${(F,\widetilde{\rho})}$ is a Morse-Smale pair on $Q$
    \item ${\widetilde{\rho}|_{L_{\pm}}=\rho_{\pm}}$
    \item ${(F,\widetilde{\rho})|_{\pi^{-1}(z)}=(f_{\pm}+c(z),\rho_{\pm})}$ near ${z= \pm 1}$
    \item ${\pi_*(\mathrm{grad}_{\widetilde{\rho}}F)}$ is constant on ${\pi^{-1}(z)}$ for any ${z\in \partial D^2}$.
\end{itemize}
We choose two flows ${\gamma_{\pm}:\mathbb{R}\rightarrow \partial D^2}$ of the vector field ${\pi_*(\mathrm{grad}_{\widetilde{\rho}}F)}$ which connect ${+1}$ and ${-1}$ so that ${\pm \mathrm{Im}(\gamma_{\pm}(t))>0}$ holds. Let ${(H_t^{\pm},h_t^{\pm})}$ be two families of pairs of functions and Riemannian metrics on $L$ defines as follows:
\begin{gather*}
        H_t^{-}=F|_{\pi^{-1}(\gamma(t))}  \\
        h_t^{-}=\widetilde{\rho}|_{\pi^{-1}(\gamma(t))}
\end{gather*}
\begin{gather*}
    H_t^{+}=\begin{cases}
        F|_{\pi^{-1}(\gamma(t))} & \mathrm{Re}\gamma(t)\le 0 \\
        (g_1)_* F|_{\pi^{-1}(\gamma(t))} & \mathrm{Re}\gamma(t)\ge 0
    \end{cases}
\end{gather*}
\begin{gather*}
    h_t^{+}=\begin{cases}
        \widetilde{\rho}|_{\pi^{-1}(\gamma(t))} & \mathrm{Re}\gamma(t)\le 0  \\
        (g_1)_*\widetilde{\rho}|_{\pi^{-1}(\gamma(t))} & \mathrm{Re}\gamma(t)\ge 0 .
    \end{cases}
\end{gather*}
Note that
\begin{gather*}
    (H_t^{-},h_t^{-})=\begin{cases}
        (f_++c(\gamma_-(t)),\rho_+)  & t \lll 0  \\
        (f_-+c(\gamma_-(t)),\rho_-)  & t \ggg 0
    \end{cases}
\end{gather*}
and 
\begin{gather*}
    (H_t^+,h_t^+)=\begin{cases}
        ((g_1)_*f_++c(\gamma_+(t)),(g_1)_*\rho_+) & t\lll 0  \\
        (f_-+c(\gamma_+(t)),\rho_-)  & t\ggg 0
    \end{cases}
\end{gather*}
hold. Let 
\begin{gather*}
    \Phi:CM_*(f_+,\rho_+)\longrightarrow CM_*(f_-,\rho_-)  \\
    x\longmapsto \sum \sharp_{\mathbb{Z}_2} \mathcal{N}(x,y:(H_t^{-},h_t^{-}))\cdot y
\end{gather*}
be the canonical continuation map. Let ${\Phi^{g_1}}$ be a chain map defined as follows:
\begin{gather*}
    \Phi^{g_1}:CM_*(f_+,\rho_+)\longrightarrow CM_*(f_-,\rho_-) \\
    x\longmapsto \sum \sharp_{\mathbb{Z}_2} \mathcal{N}(g_1(x),y:(H_t^{+},h_t^+))\cdot y .
\end{gather*}
Note that ${\Phi^{g_1}}$ is a composition of the natural $g_1$-action
\begin{gather*}
    CM(f_+,\rho_+)\longrightarrow CM_*((g_1)_*f_+,(g_1)_*\rho_+) \\
    x\longmapsto g_1(x)
\end{gather*}
and the canonical continuation map
\begin{gather*}
    CM_*((g_1)_*f_+,(g_1)_*\rho_+)\longrightarrow CM_*(f_-,\rho_-) \\
    z\longmapsto \sum \sharp_{\mathbb{Z}_2} \mathcal{N}(z,y:(H_t^+,h_t^+))\cdot y.
\end{gather*}
Morse chain complex of $F$ is a direct sum of ${CM_*(f_-,\rho_-)}$ and ${CM_*(f_+,\rho_+)}$.
\begin{gather*}
    CM_*(F,\widetilde{\rho})=CM_*(f_-,\rho_-)\oplus CM_{*-1}(f_+,\rho_+)
\end{gather*}
Morse boundary operator ${\partial ^M_{F}}$ of ${ CM_*(F,\widetilde{\rho})}$ is described as follows. We denote the boundary operator on ${CM_*(f_{\pm},\rho_{\pm})}$ by ${\partial ^M_{f_{\pm}}}$.
\begin{gather*}
    \partial^M_{F}:CM_*(f_-,\rho_-)\oplus CM_{*-1}(f_+,\rho_+)\longrightarrow CM_{*-1}(f_-,\rho_-)\oplus CM_{*-2}(f_+,\rho_+) \\
    (x_-,x_+)\longmapsto (\partial^M_{f_-}(x_-)+\Phi(x_+)+\Phi^{g_1}(x_+),\partial^M_{f_+}(x_+))
\end{gather*}
This induces the following exact sequence:
\begin{gather*}
    \cdots \longrightarrow HM_*(f_+,\rho_+)\stackrel{k}{\longrightarrow}HM_*(f_-,\rho_-)\stackrel{i}{\longrightarrow}HM_*(F,\widetilde{\rho})\\
    \stackrel{j}{\longrightarrow} HM_{*-1}(f_+,\rho_+)\longrightarrow \cdots
\end{gather*}
Maps $i$,$j$ and $k$ are induced from the following chain maps:
\begin{gather*}
    i:CM_*(f_-,\rho_-)\longrightarrow CM_*(f_-,\rho_-)\oplus CM_{*-1}(f_+,\rho_+)  \\
    x_-\longmapsto (x_-,0)
\end{gather*}
\begin{gather*}
    j:CM_*(f_-,\rho_-)\oplus CM_{*-1}(f_+,\rho_+)\longrightarrow CM_{*-1}(f_+,\rho_+) \\
    (x_-,x_+)\longmapsto x_+
\end{gather*}
\begin{gather*}
    k:CM_*(f_+,\rho_+)\longrightarrow CM_*(f_-,\rho_-) \\
    x_+\longmapsto \Phi(x_+)+\Phi^{g_1}(x_+)
\end{gather*}
The exactness of the above long exact sequence is easy to check. In particular, assume that ${f_+-f_-}$ is constant and ${\rho_+=\rho_-}$ hold. In this case we can identify ${CM_*(f_+,\rho_+)}$ and ${CM_*(f_-,\rho_-)}$ (we denote this by ${CM_*(f,\rho)}$). We get the following long exact sequence.

\[
\begin{tikzcd}
HM(f,\rho) \arrow[rr,"(g_1)_*-\mathrm{Id}"] &               & HM(f,\rho) \arrow[ldd,"i"] \\
             &               &               \\
             & HM(F,\widetilde{\rho}) \arrow[luu,"\lbrack -1\rbrack"] &              
\end{tikzcd}
\]

\begin{gather*}
     \cdots \longrightarrow HM_*(f,\rho)\stackrel{[g_1]-\mathrm{Id}}{\longrightarrow}HM_*(f,\rho)\longrightarrow HM_*(F,\widetilde{\rho})\\
    \longrightarrow HM_{*-1}(f,\rho)\longrightarrow \cdots
\end{gather*}
This is a Morse theoretical description of Wang exact sequence.

\section{Seidel-type map and proof of Theorem 1}
We fix a Hamiltonian isotopy ${\mathbf{g}=\{\phi_t\}_{t\in[0,1]}\subset \mathrm{Ham}(M,\omega)}$ so that ${\phi_0=\mathrm{Id}}$ and ${\phi_1\in \mathrm{Ham}_L(M,\omega)}$ hold. As considered in the previous section, we have the following symplectic fibration
\begin{gather*}
    \pi:P\longrightarrow D^2  \\
    P=D_-^2\times M\coprod D_+^2\times M/ \sim
\end{gather*}
and a Lagrangian fibration
\begin{gather*}
    \pi:Q\longrightarrow \partial D^2 \\
    Q=\partial^{-}D^2\times L\coprod \partial^+D^2\times L / \sim .
\end{gather*}
Let $Z_+$ and $Z_-$ be positive and negative cylinders:
\begin{gather*}
    Z_+=[0,+\infty)\times [0,1]\\
    Z_-=(-\infty,0]\times [0,1].
\end{gather*}
We identify $D_+\backslash \{1\}$ and $Z_+$, ${D_-\backslash \{-1\}}$ and ${Z_-}$. We fix a cut-off function $\beta$ as follows:

\begin{gather*}
    \beta:[0,+\infty)\longrightarrow [0,1] 
\end{gather*}
\begin{gather*}
    \beta(s)=\begin{cases}
        1 & 0\le s\le 1 \\ 0 & s\ge 2 .
    \end{cases}
\end{gather*}
Then, we can define a $2$-form ${\omega_{\mathbf{g}}}$ on $P$ as follows. Note that we choose a coordinate in the decomposition
\begin{gather*}
     P=D_-^2\times M\coprod D_+^2\times M/ \sim
\end{gather*}
and we identify $D_+\backslash \{1\}$ and $Z_+$, ${D_-\backslash \{-1\}}$ and ${Z_-}$. 

\begin{gather*}
    \omega_{\mathbf{g}}((s,t),x)=\begin{cases}
        \omega & ((s,t),x)\in Z_-\times M  \\
        \omega-\beta(s)d\overline{H}_tdt-\beta'(s)\overline{H}_tdsdt & ((s,t),x)\in Z_+\times M
    \end{cases}
\end{gather*}
${\overline{H}_t}$ is a Hamiltonian function which generate a Hamiltonian isotopy ${\{(\phi_H^t)^{-1}\}}$.
This $\omega_{\mathbf{g}}$ is a coupling form of the fibration
\begin{gather*}
    \pi:P\longrightarrow D^2.
\end{gather*}
In other words, $\omega_{\mathbf{g}}$ satisfies the following conditions:
\begin{itemize}
\item  $\omega_{\mathbf{g}}|_{\pi^{-1}(x)}=\omega \ \ \ \forall x\in D^2  $
\item  $d\omega_{\mathbf{g}}=0$ \ ($\omega_{\mathbf{g}}$ is a closed form.)
\item  $\pi_*(\omega_{\mathbf{g}}^{(n+1)})=0\in \Omega^2(D^2)$
\end{itemize}

Let ${(f_+,\rho_+)}$ be a Morse-Smale pair on $L_+$ and let ${(f_-,\rho_-)}$ be a Morse-Smale pair on $L_-$. We define the following ring $\Omega$.
\begin{gather*}
    \Omega=\Big\{ \ \sum_{\lambda\in \mathbb{R}}a_{\lambda}T^{\lambda} \ \Big| \ a_{\lambda}\in \mathbb{Z}_2 \ \Big\}
\end{gather*}
Unlike the usual Novikov ring, an element of $\Omega$ is a finite sum. We consider $\Omega$-coefficient Morse chain complex and $\Omega$-coefficient Mose homology.
\begin{gather*}
    CM(f_{\pm},\rho_{\pm}:\Omega)= CM(f_{\pm},\rho_{\pm})\otimes_{\mathbb{Z}_2} \Omega 
\end{gather*}
\begin{gather*}
    HM(f_{\pm},\rho_{\pm}:\Omega)=H(CM(f_{\pm},\rho_{\pm}:\Omega),\partial)
\end{gather*}
The first thing to do is to use fibrations
\begin{gather*}
    \pi:(P,Q)\longrightarrow (D^2,\partial D^2)
\end{gather*}
to construct a Seidel-type map
\begin{gather*}
    \Phi_{\mathbf{g}}:HM((f_-,\rho_-):\Omega)\longrightarrow HM((f_+,\rho_+):\Omega).
\end{gather*}

\vspace{7mm}
\begin{Lem}
    For any ${\epsilon>0}$, there is a positive function
    \begin{gather*}
        \kappa: D^2\longrightarrow \mathbb{R}_{> 0}
    \end{gather*}
    such that ${\pi^*(\kappa \ dxdy)+\omega_{\mathbf{g}}}$ is a symplectic form on $P$ and 
\begin{gather*}
    \int_{D^2}\kappa(x,y) dxdy \le ||H||_{+}+\epsilon
\end{gather*}
holds.
\end{Lem}
\vspace{7mm}

\textbf{Proof}  \ Recall that ${\omega_{\mathbf{g}}}$ is trivial over the half disc $D_-^2$:
\begin{gather*}
    \omega_{\mathbf{g}}|_{\pi^{-1}(D_-^2)}=\omega .
\end{gather*}
So it suffices to prove the lemma over ${D_+^2}$ (or ${Z_+\cong D_+^2\backslash \{1\}}$). The following equalities hold on ${\pi^{-1}(D_+^2\backslash \{1\})\cong Z_+\times M}$. Let ${\overline{\kappa}}$ be a positive function on $Z_+$.

\begin{gather*}
    \iota_{\frac{\partial}{\partial s}}(\overline{\kappa}(s,t)dsdt+\omega_{\mathbf{g}})=\overline{\kappa}(s,t)dt-\beta'(s)\overline{H}dt
\end{gather*}
\begin{gather*}
    \iota_{\frac{\partial}{\partial t}}(\overline{\kappa}(s,t)dsdt+\omega_{\mathbf{g}})=-\overline{\kappa}(s,t)+\beta(s) d\overline{H}dt+\beta'(s)\overline{H}ds
\end{gather*}
These equalities imply that ${\overline{\kappa}dsdt+\omega_{\mathbf{g}}}$ is non-degenerate (hence symplectic) if 
\begin{gather*}
    \overline{\kappa}(s,t)>|\beta'(s)|\max H_t
\end{gather*}
holds. We choose ${\overline{\kappa}:Z_+\rightarrow \mathbb{R}_{> 0}}$ so that 
\begin{gather*}
    \int_{Z_+}\overline{\kappa}(s,t)dsdt\le ||H||_++\frac{1}{2}\epsilon \\
    \overline{\kappa}(s,t)>|\beta'(s)|\max H_t
\end{gather*}
holds. Then we can choose ${\kappa:D^2\rightarrow \mathbb{R}_{\ge 0}}$ so that
\begin{gather*}
    \kappa dxdy|_{D_+^2}=\overline{\kappa}dsdt  \\
    \int_{D^2}\kappa dxdy\le ||H||_++\epsilon
\end{gather*}
holds. Then ${\pi^*(\kappa dxdy)+\omega_{\mathbf{g}}}$ is a symplectic form on $P$.
\begin{flushright}
    $\square$
\end{flushright}
From now on, we fix a sufficiently small ${\epsilon>0}$ and a positive function ${\kappa:D^2\rightarrow \mathbb{R}_{>0}}$ as in Lemma 1. We define the space of almost complex structures on P.
\begin{gather*}
    \mathcal{J}(P,\kappa)=\Bigg\{ \ J\in \mathrm{End}(TP) \ \Bigg| \ \begin{matrix}
        J^2=-\mathrm{Id}, \ d\pi\circ J=\sqrt{-1}\circ d\pi  \\ J:\mathrm{compatible \ with \ }\pi^*(\kappa dxdy)+\omega_{\mathbf{g}}
    \end{matrix} \ \Bigg\}
\end{gather*}
We define ${\widetilde{\pi}_2(P,Q)}$ as follows. Let ${[D^2]\in H_2(D^2,\partial D^2)}$ be the relative fundamental class.
\begin{gather*}
    \widetilde{\pi}_2(P,Q)=\{u:(D^2,\partial D^2)\rightarrow (P,Q) \ | \ \pi_*u_*([D^2])=[D^2]\}/\sim
\end{gather*}
The equvalence relation ${\sim}$ is defined by homotopies. $u\sim v$ holds if and only if there exists a $1$-parameter family 
\begin{gather*}
    u_t:(D^2,\partial D^2)\longrightarrow (P,Q) \ \ \ (t\in [0,1])
\end{gather*}
so that ${u_0=u}$ and ${u_1=v}$ hold. Note that 
\begin{gather*}
    \langle \omega_{\mathbf{g}},[u]\rangle=\int_{D^2}u^*\omega_{\mathbf{g}}
\end{gather*}
is well-defined for ${[u]\in \widetilde{\pi}_2(P,Q)}$ because $Q$ is a Lagrangian submanifold of $P$. We consider the followig moduli space of J-holomorphic sections for ${J\in \mathcal{J}(P,\kappa)}$ and ${\alpha \in \widetilde{\pi}_2(P,Q)}$.

\begin{gather*}
    \widetilde{\mathcal{M}}(J,\alpha)=\Bigg\{ \ v:(D^2,\partial D^2)\rightarrow (P,Q) \ \Bigg| \ \begin{matrix}
       J\circ dv=dv  \circ \sqrt{-1}  \\ [v]=\alpha
    \end{matrix} \ \Bigg\}
\end{gather*}
The automorphism group of $D^2$ which preserves ${\pm 1}$ is isomorphic to $\mathbb{R}$. When we identify ${D^2\backslash \{\pm 1\}}$ to ${\mathbb{R}\times [0,1]}$, this is the $\mathbb{R}$-shift of the first factor. We divide ${\widetilde{\mathcal{M}}(J,\alpha)}$ by this $\mathbb{R}$-action.
\begin{gather*}
    \mathcal{M}(J,\alpha)=\widetilde{\mathcal{M}}(J,\alpha)/\mathbb{R}
\end{gather*}

\vspace{7mm}
\begin{Lem}
    If $\mathcal{M}(J,\alpha)$ is not empty, the following inequality is satisfied.
    \begin{gather*}
        \langle \omega_{\mathbf{g}}, \alpha \rangle \ge -||H||_+-\epsilon
    \end{gather*}
\end{Lem}
\vspace{7mm}
\textbf{Proof} \ We fix ${u\in \widetilde{\mathcal{M}}(J,\alpha)}$. Recall that $J$ is ${\pi^*(\kappa dxdy)+\omega_{\mathbf{g}}}$ compatible almost complex structure on $P$. So,
\begin{gather*}
    \int_{D^2}u^*\{\pi^*(\kappa dxdy)+\omega_{\mathbf{g}}\}\ge 0
\end{gather*}
is satisfied. This implies that the following inequality holds.
\begin{gather*}
    \int_{D^2}u^*\omega_{\mathbf{g}}\ge -\int_{D^2}\kappa(x,y)dxdy\ge -(||H||_++\epsilon)
\end{gather*}
\begin{flushright}
    $\square$
\end{flushright}

For ${x_-\in \mathrm{Crit}(f_-)}$ and ${x_+\in \mathrm{Crit}(f_+)}$, we define the following moduli space:
\begin{gather*}
    \mathcal{N}(x_-,x_+:J,\alpha)=\Bigg\{ \ u\in \mathcal{M}(J,\alpha) \ \Bigg| \ \begin{matrix}u(-1)\in W^u(x_-)    \\ u(+1)\in W^s(x_+)\end{matrix} \ \Bigg\}.
\end{gather*}
We use this moduli space to define a Seidel-type chain map $\Phi_{\mathbf{g}}$:

\begin{gather*}
    \Phi_{\mathbf{g}}:CM(f_-,\rho_-:\Omega)\longrightarrow CM(f_+,\rho_+:\Omega)  
\end{gather*}
\begin{gather*}
   \Phi_{\mathbf{g}} (x_-)\stackrel{\mathrm{def}}{=} \\  \sum_{\langle \omega_{\mathbf{g}},\alpha \rangle \le ||H||_-+\epsilon}\sum_{x_+\in \mathrm{Crit}(f_+)}\sharp_{\mathbb{Z}_2} \mathcal{N}(x_-,x_+:J,\alpha)T^{-\langle \omega_{\mathbf{g},\alpha \rangle}}\cdot x_+ .
\end{gather*}
This induces a map between homologies.
\begin{gather*}
     \Phi_{\mathbf{g}}:HM(f_-,\rho_-:\Omega)\longrightarrow HM(f_+,\rho_+:\Omega)  
\end{gather*}
Next, we consicer a Seidel-type map for the inverse isotopy 
\begin{gather*}
\mathbf{g}^{-1}=\{\phi_{\overline{H}}^t\}\subset \mathrm{Ham}(M,\omega).
\end{gather*}
We have a symplectic fibration and a Lagrangian fibration as before.
\begin{gather*}
    \pi:(P_{\mathbf{g}^{-1}},Q_{\mathbf{g}^{-1}})\longrightarrow (D^2,\partial D^2)
\end{gather*}
As in the case of $\mathbf{g}$, we can define a Seidel-type map $\Phi_{\mathbf{g}^{-1}}$ for $\mathbf{g}^{-1}$. We choose ${\kappa':D^2\rightarrow \mathbb{R}_{>0}}$ instead of ${\kappa}$ so that
\begin{gather*}
    \int_{D^2}\kappa'dxdy\le ||H||_-+\epsilon
\end{gather*}
holds and ${\pi^*(\kappa'dxdy)+\omega_{\mathbf{g}^{-1}}}$ is a symplectic form on ${P_{\mathbf{g}^{-1}}}$. We define ${\mathcal{J}'(P_{\mathbf{g}^{-1}},\kappa')}$, ${\mathcal{M}'(J',\beta)}$ and ${\mathcal{N}'(x_+,x_-:J',\beta)}$ as in the case of ${\mathbf{g}}$. Then, ${\Phi_{\mathbf{g}^{-1}}}$ can be written in the following form.
\begin{gather*}
    \Phi_{\mathbf{g}^{-1}}:CM(f_+.\rho_+:\Omega)\longrightarrow CM(f_-,\rho_-:\Omega)
\end{gather*}
\begin{gather*}
    \Phi_{\mathbf{g}^{-1}}(x_+)= \\
    \sum_{\langle \omega_{\mathbf{g}^{-1}}, \beta \rangle\le ||H||_++\epsilon}\sum_{x_-\in \mathrm{Crit}(f_-)}\sharp_{\mathbb{Z}_2} \mathcal{N}'(x_+,x_-:J',\beta)T^{-\langle \omega_{\mathbf{g}^{-1}},\beta \rangle}\cdot x_-
\end{gather*}
We study the difference of ${\Phi_{\mathbf{g}^{-1}}\circ \Phi_{\mathbf{g}}}$ and ${\mathrm{Id}}$. We apply chain homotopies constructed by standard glueing arguments in Morse theory and Floer theory. We define ${\mathrm{ev}_0}$ and ${\mathrm{ev}_1}$ as follows:
\begin{gather*}
    \mathrm{ev}_0:\mathcal{M}(J,\alpha)\longrightarrow L_+  \\
    [u]\longmapsto u(+1)
\end{gather*}
\begin{gather*}
    \mathrm{ev}_1:\mathcal{M}'(J',\alpha)\longrightarrow L_+  \\
    [v]\longmapsto v(-1) .
\end{gather*}
We define the fiber product of ${\mathcal{M}(J,\alpha)}$ and ${\mathcal{M}'(J',\beta)}$:
\begin{gather*}
    \mathcal{M}(J,\alpha)\times_{L_+}\mathcal{M}'(J',\beta)= \\ \Big\{ \ (u,v)\in \mathcal{M}(J,\alpha)\times \mathcal{M}'(J',\beta) \ \Big| \ \mathrm{ev}_0(u)=\mathrm{ev}_1(v) \ \Big\}.
\end{gather*}
For any ${y,z\in \mathrm{Crit}(f_-)}$, we define the following moduli space:
\begin{gather*}
    \mathcal{N}(y,z:\alpha,\beta,J,J')=\\ \bigg\{\ (u,v)\in  \mathcal{M}(J,\alpha)\times_{L_+}\mathcal{M}'(J',\beta) \ \bigg| \ \begin{matrix}
        u(-1)\in W^u(y)  \\ v(+1)\in W^s(z) 
    \end{matrix}\ \bigg\}.
\end{gather*}
We define another chain map $\Psi$ as follows:
\begin{gather*}
    \Psi:CM(f_-,\rho_-:\Omega)\longrightarrow CM(f_-,\rho_-:\Omega)
\end{gather*}
\begin{gather*}
    \Psi(y)=\\ \sum_{\langle \omega_{\mathbf{g}},\alpha \rangle\le ||H||_-+\epsilon}   \sum_{\langle \omega_{\mathbf{g}^{-1}},\beta \rangle \le ||H||_++\epsilon}\sum_{z\in\mathrm{Crit}(f_-)}\sharp_{\mathbb{Z}_2} \mathcal{N}(u,z:\alpha,\beta,J,J')T^{-\langle \omega_{\mathbf{g}},\alpha \rangle -\langle \omega_{\mathbf{g}^{-1}},\beta \rangle}\cdot z .
\end{gather*}
First, we study the difference of $\Psi$ and ${\mathrm{Id}}$. By applying gluing of J-homolorphic discs in ${\mathcal{N}(y,z:\alpha,\beta,J,J')}$, we can get a J-holomorphic disc. Let ${J_R}$ be an almost complex structure on ${D^2\times M}$ for ${R\gg 0}$ as follows. We identify ${D^2\backslash \{\pm 1\}}$ and ${\mathbb{R}\times [0,1]}$.
\begin{gather*}
    J_R((s,t),x)=\begin{cases}
        J((s+R,t),x)  & s\le 0 \\ J((s-R,t),x)  & s\ge 0
    \end{cases}
\end{gather*}
Then, gluing of J-holomorphic discs (see Theorem 4.1.2 in \cite{BC}) implies that for ${(u,v)\in \mathcal{N}(y,z:\alpha,\beta,J,J')}$, ${\langle \omega_{\mathbf{g}},\alpha \rangle\le ||H||_++\epsilon}$, ${\langle \omega_{\mathbf{g}^{-1}},\beta \rangle \le ||H||_-+\epsilon}$, there is a $1$-parameter family of ${u_R:\mathbb{R}\times [0,1]\rightarrow (M,L)}$ as follows:
\begin{gather*}
\begin{cases}
    \partial_su_R(s,t)+J_R(\partial_tu_R(s,t))=0  \\
    u_R \mathrm{\ converves \ to \ }(u,v) \ \ \ (u\rightarrow +\infty)
\end{cases}.
\end{gather*}
In other words, $u_R$ is obtained by gluing ${u\in \mathcal{M}(J,\alpha)}$ and ${v\in \mathcal{M}'(J',\beta)}$ along ${u(+1)=v(-1)\in L_+}$. For ${y,z\in \mathrm{Crit}(f_-)}$, ${\gamma\in \pi_2(M,L)}$ and ${R\gg 0}$, we define ${\mathcal{S}(y,z:\gamma,R)}$ as follows: 
\begin{gather*}
    \mathcal{S}(y,z:\gamma,R)= \\
    \bigg\{\ u:(D^2,\partial D^2)\rightarrow(M,L) \ \bigg| \ \begin{matrix}
        [u]=\gamma, \ \partial _su+J_R\partial_tu=0 \\ u(-1)\in W^u(y), \ u(+1)\in W^s(z)
    \end{matrix} \bigg\} .
\end{gather*}

\vspace{7mm}
\begin{Rem}
We have two cupling forms $\omega_{\mathbf{g}}$ on $P$ and $\omega_{\mathbf{g}^{-1}}$ on ${P_{\mathbf{g}^{-1}}}$. By identifying ${P\sharp P_{\mathbf{g}^{-1}}}$ and ${D^2\times M}$, we get a coupling form ${\omega_{\mathbf{g}}\sharp \omega_{\mathbf{g}^{-1}}}$. Let ${D_-^2}$ be the unit disc with opposite orientation. By gluing ${D^2\times M}$ and ${D_-^2\times M}$ along the boundary, we get a trivial symplectic fibration ${S^2\times M\rightarrow M}$ and a coupling form ${\tau}$ on ${S^2\times M}$ defined by
\begin{gather*}
    \tau=\begin{cases}
        \omega & \mathrm{on} \ D^2\times M \\ \omega_{\mathbf{g}}\sharp \omega_{\mathbf{g}^{-1}} & \mathrm{on} \ D_-^2\times M
    \end{cases}.
\end{gather*}
Then, the uniqueness of the coupling class implies that ${[\tau]=[\omega]}$ holds on ${S^2\times M}$. For any section ${s:D^2\rightarrow D^2\times M}$, we get a section
\begin{gather*}
    s\sharp \overline{s}:S^2\longrightarrow S^2\times M
\end{gather*}
by gluing two $s$ along the boundary. Then we see that
\begin{gather*}
    \int_{D^2}s^*\omega-\int_{D^2}s^*(\omega_{\mathbf{g}}\sharp \omega_{\mathbf{g}^{-1}})=
    \int_{S^2}(s\sharp \overline{s})^*\tau =  
        \int_{S^2}(s\sharp \overline{s})^*\omega=0
\end{gather*}
holds. In particular, we see that 
\begin{gather*}
    \int_{D^2}s^*(\omega_{\mathbf{g}}\sharp \omega_{\mathbf{g}^{-1}})=\int_{D^2}s^*\omega
\end{gather*}
hodls for any section ${s:D^2\rightarrow D^2\times M}$.
\end{Rem}
\vspace{7mm}
We apply ${\mathcal{S}(y,z:\gamma,R)}$ to study the difference ${\Psi-\mathrm{Id}}$. We consider a ${1}$-parameter family of ${\mathcal{S}(y,z:\gamma,R)}$ parametrized by ${r\ge R\gg 0}$ as follows:
\begin{gather*}
    \mathcal{S}^{\ge R}(y,z:\gamma)=\coprod_{r\ge R}\mathcal{S}(y,z:\gamma,r).
\end{gather*}
Then we can define the following two maps $\mathrm{I}_R$ and ${h}$.
\begin{gather*}
    \mathrm{I}_R:CM(f_-,\rho_-:\Omega)\longrightarrow CM(f_-,\rho_-:\Omega) \\
    y\longmapsto \sum_{\langle \omega, \gamma \rangle \le||H||+2\epsilon}\sum_{z\in\mathrm{Crit}(f_-)}\sharp_{\mathbb{Z}_2} \mathcal{S}(y,z:\gamma,R)T^{-\langle \omega,\gamma\rangle}\cdot z
\end{gather*}
\begin{gather*}
    h:CM(f_-,\rho_-:\Omega)\longrightarrow CM(f_-,\rho_-:\Omega) \\
    y\longmapsto \sum_{\langle \omega, \gamma \rangle \le||H||+2\epsilon}\sum_{z\in\mathrm{Crit}(f_-)}\sharp _{\mathbb{Z}_2}\mathcal{S}^{\ge R}(y,z:\gamma)T^{-\langle \omega,\gamma\rangle}\cdot z
\end{gather*}
We denote the followig projection to the non-negative part and non-positive part by ${\Pi_+}$ and ${\Pi_-}$.

\begin{gather*}
    \Pi_+:CM(f_-,\rho_-:\Omega)\longrightarrow CM(f_-,\rho_-:\Omega)  \\
    \sum_{z\in \mathrm{Crit}(f_-)}\sum_{\lambda\in \mathbb{R}}a_{\lambda,z}T^{\lambda}\cdot z \longmapsto  \sum_{z\in \mathrm{Crit}(f_-)}\sum_{\lambda\ge 0}a_{\lambda,z}T^{\lambda}\cdot z 
\end{gather*}

\begin{gather*}
    \Pi_-:CM(f_-,\rho_-:\Omega)\longrightarrow CM(f_-,\rho_-:\Omega)  \\
    \sum_{z\in \mathrm{Crit}(f_-)}\sum_{\lambda\in \mathbb{R}}a_{\lambda,z}T^{\lambda}\cdot z \longmapsto  \sum_{z\in \mathrm{Crit}(f_-)}\sum_{\lambda\le 0}a_{\lambda,z}T^{\lambda}\cdot z 
\end{gather*}
Next lemma is inspired by Chekanov's paper \cite{C}, where Chekanov proved that Lagrangian intersection property of ``small" Hamiltonian isopoties.
\vspace{7mm}
\begin{Lem}
    The following equality holds.
    \begin{gather*}
        \Pi_+(\mathrm{I}_R+\Psi+h\partial +\partial h)\Pi_-=0
    \end{gather*}
\end{Lem}
\vspace{7mm}
\textbf{Proof} \ The proof relies on the following standard cobordism argument in Floer and Morse theory. We study the boudary of ${\mathcal{S}^{\ge R}(y,z:\gamma)}$ when ${\langle \omega,\gamma \rangle \le 0}$ holds and dimension of ${\mathcal{S}^{\ge R}(y,z:\gamma)}$ is equal to $1$. Let ${u_r}$ be a $1$-parameter family of solutions of the equation
\begin{gather*}
    u_r:(D^2,\partial D^2)\longrightarrow (M,L)  \\
    \langle \omega, u_r\rangle \le 0 \\
    \partial_su_r(s,t)+J_r(\partial_tu_r(s,t))=0
\end{gather*}
so that when ${r\rightarrow +\infty}$, ${u_r}$ splits into ${u\in \mathcal{M}(J,\alpha)}$ and ${v\in \mathcal{M}'(J',\beta)}$. We prove that ${\langle\omega_{\mathbf{g}},\alpha\rangle\le ||H||_-+\epsilon}$ and ${\langle \omega_{\mathbf{g}^{-1}},\beta \rangle\le ||H||_++\epsilon}$ hold. Assume that ${\langle\omega_{\mathbf{g}},\alpha\rangle> ||H||_-+\epsilon}$ holds. In this case,
\begin{gather*}
    \langle \omega_{\mathbf{g}^{-1}},\beta \rangle=\langle \omega,u_r\rangle-\langle \omega_{\mathbf{g}^{-1}},\beta \rangle<-||H||_--\epsilon
\end{gather*}
holds and this contradicts to Lemma 2 (we replace ${\mathbf{g}}$ and ${\mathbf{g}^{-1}}$ in Lemma 2). Similarly, we see that ${\langle \omega_{\mathbf{g}^{-1}},\beta \rangle\le ||H||_++\epsilon}$ holds. This implies that the boundary of ${\mathcal{S}^{\ge R}(y,z:\gamma)}$ (when ${\langle \omega,\gamma\rangle \le 0}$ holds and the dimension of ${\mathcal{S}^{\ge R}(y,z:\gamma)}$ is equal to $1$) can be written as follows:
\begin{gather*}
    \partial \mathcal{S}^{\ge R}(y,z:\gamma)= \\
    \coprod_{y'\in\mathrm{Crit}(f_-)}\mathcal{M}(y,y':(f_-,\rho_-))\times \mathcal{S}^{\ge R}(y',y:\gamma) \\
    +\coprod_{z'\in\mathrm{Crit}(f_-)}\mathcal{S}^{\ge R}(y,z':\gamma)\times \mathcal{M}(z',z:(f_-,\rho_-)) \\
    +\mathcal{S}(y,z:\gamma,R)  \\
    +\coprod_{\alpha\sharp \beta=\gamma}\coprod_{\langle\omega_{\mathbf{g}},\alpha\rangle\le ||H||_-+\epsilon}\coprod_{\langle \omega_{\mathbf{g}^{-1}},\beta \rangle\le ||H||_++\epsilon}\mathcal{N}(y,z:\alpha,\beta,J,J') .
\end{gather*}
This decomposition implies that 
\begin{gather*}
     \Pi_+(\mathrm{I}_R+\Psi+h\partial +\partial h)\Pi_-=0
\end{gather*}
holds.
\begin{flushright}
    $\square$
\end{flushright}
Note that $\mathrm{I}_R$ and ${\mathrm{Id}}$ are chain homotopic (it suffices to deform ${J_R}$ for constant almost complex structure $J$) and $\Psi$ and ${\Phi_{\mathbf{g}^{-1}}\Phi_{\mathbf{g}}}$ are chain homotopic. So, we proved the following lemma.
\vspace{7mm}
\begin{Lem}
    There is a map
    \begin{gather*}
        H:CM(f_-,\rho_-:\Omega)\longrightarrow CM(f_-,\rho_-:\Omega)
    \end{gather*}
    such that
    \begin{gather*}
        \Pi_+(\mathrm{Id}+\Phi_{\mathbf{g}^{-1}}\Phi_{\mathbf{g}}+H\partial+\partial H)\Pi_-=0 
    \end{gather*}
    holds.
\end{Lem}
\vspace{7mm}
\begin{Cor}
A chain map 
\begin{gather*}
    \Phi_{\mathbf{g}}^0:CM(f_-,\rho_-)\longrightarrow CM(f_+,\rho_+:\Omega) \\
     x\longmapsto \Phi_{\mathbf{g}}(x\otimes 1)
\end{gather*}
induces an injection 
\begin{gather*}
    \Phi_{\mathbf{g}}^0:HM(f_-,\rho_-)\longrightarrow HM(f_+,\rho_+:\Omega) .
\end{gather*}
\end{Cor}
\vspace{7mm}
\textbf{Proof} \ Assume that ${c\in CM(f_-,\rho_-)}$ is a cycle so that ${[c]\neq 0}$ and ${\Phi_{\mathbf{g}}^0([c])=0}$ hold. This implies that ${\Phi_{\mathbf{g}^{-1}}\Phi_{\mathbf{g}}(c\otimes 1)}$ is a boundary element in ${CM(f_-,\rho_-:\Omega)}$. Lemma 4 implies that 
\begin{gather*}
    c\otimes 1=\Pi_+(\Phi_{\mathbf{g}^{-1}}\Phi_{\mathbf{g}}(c\otimes 1)+\partial H(c\otimes 1))
\end{gather*}
and hence $c$ is a boundary element. This is a contradiction.
\begin{flushright}
    $\square$
\end{flushright}
Next, we prove Theorem 1. It suffices to prove that
\begin{gather*}
    i:HM(f_-,\rho_-)\longrightarrow HM(F,\widetilde{\rho})
\end{gather*}
is injective. Here, ${(F,\widetilde{\rho})}$ is a Morse-Smale pair on $Q$ constructed in the previous section. We extend 
\begin{gather*}
    \Phi_{\mathbf{g}}:HM(f_-,\rho_-:\Omega)\longrightarrow HM(f_+,\rho_+:\Omega)
\end{gather*}
to ${HM(F,\widetilde{\rho})}$. For ${x\in\mathrm{Crit}(F)}$, ${y_+\in\mathrm{Crit}(f_+)}$, ${\alpha\in \widetilde{\pi}_2(P,Q)}$, we consider the following moduli space:
\begin{gather*}
    \widetilde{\mathcal{N}}(x,y_+:\alpha,J)=\bigg\{\ u\in \mathcal{M}(\alpha,J) \ \bigg| \ \begin{matrix}
        u(-1)\in W^u(x)  \\ u(+1)\in W^s(y_+)
    \end{matrix} \ \bigg\}.
\end{gather*}
Then we define a chain map $\widetilde{\Phi}_{\mathbf{g}}$ as follows:
\begin{gather*}
    \widetilde{\Phi}_{\mathbf{g}}:CM(F,\widetilde{\rho}:\Omega)\longrightarrow CM(f_+,\rho_+:\Omega) \\
    x\longmapsto \sum_{y_+\in \mathrm{Crit}(f_+)}\sum_{\langle\omega_{\mathbf{g}},\alpha\rangle\le ||H||_++\epsilon}\sharp_{\mathbb{Z}_2}\widetilde{\mathcal{N}}(x,y_+:\alpha,J)T^{-\langle\omega_{\mathbf{g}},\alpha\rangle}\cdot y_+ .
\end{gather*}
Note that ${\widetilde{\Phi}_{\mathbf{g}}\circ i=\Phi_{\mathbf{g}}}$ holds by definition. Assume that ${c\in CM(f_-,\rho_-)}$ is a cycle so that ${i(c)=\partial d}$ holds for some ${d\in CM(F,\widetilde{\rho})}$. Then,
\begin{gather*}
    \Phi_{\mathbf{g}}^0([c])=\Phi_{\mathbf{g}}([c\otimes 1])=\widetilde{\Phi}_{\mathbf{g}}(i([c\otimes 1]))=\widetilde{\Phi}_{\mathbf{g}}([\partial d\otimes 1])=[\partial(\widetilde{\Phi}_{\mathbf{g}}(d\otimes 1))]=0
\end{gather*}
holds. Corollary 2 implies that ${[c]=0}$ holds. In particular, $i$ is an injection and ${(\phi_1)_*=\mathrm{Id}}$ holds.
\begin{flushright}
    $\square$
\end{flushright}

\section*{Acknowledgements}
  The author gratefully acknowledges his teacher Kaoru Ono for continuous support.

\end{document}